\documentclass[amssymb, 11pt]{amsart}
\usepackage{amsmath}
\usepackage{amsthm}
\usepackage{amsfonts}
\newcommand{\C}{\mathcal{C}}
\newcommand{\yy}{\mathbb{Y}}

\newlength{\standardunitlength}
\newtheorem{prop}{Proposition}[section]

\newtheorem{lemma}[prop]{Lemma}

\newtheorem{theorem}[prop]{Theorem}

\newcommand{\ee}{\mathbb{E}}
\newcommand{\pp}{\mathbb{P}}

\begin{document}

\begin{center}
\title [Random walks on irreducible representations] {Convergence
rates of random walk on irreducible representations of finite groups}
\end{center}

\author{Jason Fulman}
\address{University of Southern California\\ Los
Angeles, CA 90089-2532}
\email{fulman@usc.edu}



\thanks{The author was partially supported by NSA grant H98230-05-1-0031 and NSF grant DMS-050391.}

\thanks{First Version: July 15, 2006. This Version: October 17, 2006.}

\begin{abstract} Random walk on the set of irreducible representations
of a finite group is investigated. For the symmetric and general
linear groups, a sharp convergence rate bound is obtained and a cutoff
phenomenon is proved. As related results, an asymptotic description of
Plancherel measure of the finite general linear groups is given, and a
connection of these random walks with quantum computing is
noted. \end{abstract}

\maketitle

\section{Introduction}

The study of convergence rates of random walk on a finite group is a
rich subject; three excellent surveys are \cite{Al}, \cite{Sal} and
\cite{D1}. A crucial role is played by random walks where the
generating measure is constant on conjugacy classes. An exact
diagonalization of such walks is possible in terms of representation
theory; this sometimes leads to sharp convergence rate bounds and even
a proof of the cut-off phenomenon. Two important cases where this has
been carried out are the random transposition walk on the symmetric
group $S_n$ \cite{DSh} and the random transvection walk on the finite
special linear group $SL(n,q)$ \cite{H}. These ``exactly solved''
Markov chains also serve as useful base chains for which the
comparison theory of \cite{DSa} can be used to analyze many other
random walks.

The current paper investigates a dual question, namely convergence
rates of random walk on $Irr(G)$, the set of irreducible
representations of a finite group $G$. Here the stationary
distribution is not the uniform distribution, but the Plancherel
measure $\pi$ on $Irr(G)$, which assigns a representation $\lambda$
probability $\frac{d_{\lambda}^2}{|G|}$, where $d_{\lambda}$ denotes
the dimension of $\lambda$. Letting $\eta$ be a (not necessarily
irreducible) representation of $G$ whose character is real valued, one
can define a Markov chain on $Irr(G)$ as follows.  From an irreducible
representation $\lambda$, one transitions to the irreducible
representation $\rho$ with probability
\[ \frac{d_{\rho} m_{\rho}(\lambda \otimes \eta)}{ d_{\lambda} d_{\eta}}\]
where $m_{\rho}(\lambda \otimes \eta)$ denotes the multiplicity of
$\rho$ in the tensor product (also called the Kronecker product) of
$\lambda$ and $\eta$. Letting $\chi$ denote the character of a
representation, the formula \[ m_{\rho}(\lambda \otimes \eta) =
\frac{1}{|G|} \sum_{g \in G} \chi^{\rho}(g) \chi^{\eta}(g)
\overline{\chi^{\lambda}(g)} \] immediately implies that this Markov
chain is reversible with respect to the Plancherel measure $\pi$.

Since these Markov chains are almost completely unexplored, we give a
detailed list of motivation for why they are worth studying:

\vspace{.2in}

(1): The same transition mechanism on $Irr(G)$ has been studied in the
closely related case where $G$ is a compact Lie group or a Lie
algebra, instead of a finite group. Then the state space $Irr(G)$ is
infinite so the questions are of a different nature than those in the
current paper. If $G=SU(2)$, then $Irr(G)$ is equal to the integers,
and the paper \cite{ER} studied asymptotics of n-step transition
probabilities. The paper \cite{BBO} used random walk on $Irr(G)$ in
the Lie algebra case, together with the Littlemann path model, to
construct Brownian motion on a Weyl chamber.

(2): Decomposing the tensor product of two elements of $Irr(G)$ (which
is what a step the random walk on $Irr(G)$ does) is just as natural as
decomposing the product of two conjugacy classes of $G$
(which is what a step in usual random walk on $G$ does). In fact for
abelian groups, $Irr(G)$ is isomorphic to $G$, and the two theories
are the same. For example the usual nearest neighbor walk on the
circle ($G=\mathbb{Z}_n$) is obtained by letting $\eta$ be the
  average of the representation closest to the trivial representation
  and its inverse, so that $\chi^{\eta}(j) = cos \left( \frac{2 \pi
    j}{n} \right)$.

Another remark which supports the idea of thinking of random walk on
$G$ and random walk on $Irr(G)$ as ``dual'' is the following. Whereas
the eigenvalues of random walk on $G$ generated by a conjugacy classes
$C$ are $\frac{\chi^{\lambda}(C)}{d_{\lambda}}$ where $C$ is fixed and
$\lambda$ varies \cite{DSh}, the eigenvalues of random walk on
$Irr(G)$ determined by the representation $\eta$ are
$\frac{\chi^{\eta}(C)}{d_{\eta}}$ where $\eta$ is fixed and $C$ varies
\cite{F2}.

(3): Decomposing the tensor product of irreducible representations of
a finite group $G$ is useful in quantum computing: see for instance
\cite{K} or page 7 of \cite{MR}, where the transition mechanism on
$Irr(G)$ is called the natural distribution of $\rho$ in $\lambda \otimes
\eta$. There are also interesting connections to free probability
theory \cite{B1}.

(4): Combinatorialists have studied the decomposition of the r-fold
tensor product of a fixed element of $Irr(G)$, when $G$ is a finite
group. Some recent results appear in \cite{GC} and \cite{GK}. Our
previous papers \cite{F1}, \cite{F3} used convergence rates of random
walk on $Irr(G)$ to study the decomposition of tensor products. The
current paper gives sharp convergence rate bounds in some cases and so
better results. It should also be noted that if $G$ is a Lie algebra,
there has been nice work done on the decomposition of the r-fold
tensor product of a fixed element of $Irr(G)$ (\cite{B2}, \cite{GM},
\cite{TZ}).

(5): There are dozens of papers written about the Plancherel measure
of the symmetric group (see the seminal papers \cite{J}, \cite{O},
\cite{BOO} and the references therein). There are many reasons why to
study a probability distribution $\pi$ of interest, it is useful to
have a Markov chain which is reversible with respect to $\pi$, and
which can be completely analyzed. First, Stein's method can then be
used to study $\pi$ (see \cite{F1},\cite{F4} where random walk on
$Irr(G)$ is used to obtain new results on Plancherel measure). Second,
convergence rate results of the chain can be used to prove
concentration inequalities for the measure $\pi$ \cite{C}. Third, the
chain can be used for Monte Carlo simulations from $\pi$, and
convergence rate results let one know how accurate the simulations
are. Fourth, the comparison theory of \cite{DSa} can be applied to
analyze the convergence rate of other chains with stationary
distribution $\pi$.

(6): The Markov chains on $Irr(G)$ are a tractable testing ground for
results in finite Markov chain theory. Whereas random walks on groups
(such as the random transposition and random transvection walk) have
simple transition probabilities but a complicated spectrum, the walks
analyzed in this paper (which are dual to the random transposition and
transvection walks) have complicated transition probabilities but a
simple spectrum. This has the effect of making convergence rate upper
bounds somewhat easier to prove, but convergence rate lower bounds
somewhat harder to prove.

\vspace{.2in}

The organization of this paper is as follows. Section \ref{prelim}
gives the needed background on Markov chain theory and defines the
cutoff phenomenon. Section \ref{groupsprelim} recalls the
diagonalization of the Markov chains on $Irr(G)$ and develops a group
theoretic tool useful for proving lower bounds. Section \ref{symirrep}
proves sharp convergence rate results for random walk on $Irr(G)$ when
$G$ is the symmetric group and $\eta$ is the defining representation
(whose character on a permutation is the number of fixed
points). Section \ref{glirrep} proves sharp convergence rate results
for random walk on $Irr(G)$ when $G$ is the finite general linear
group and $\eta$ is the permutation representation on the natural
vector space (whose character is the number of fixed vectors). Section
\ref{plancherel} obtains an elegant asymptotic description of the
Plancherel measure of $GL(n,q)$ (the next case to be understood after
the well studied case of the symmetric groups). Also, a connection
with the proof of the convergence rate lower bound of Section
\ref{glirrep} is noted. The paper closes with the very brief Section
\ref{quantum}, which uses random walk on $Irr(G)$ to explain (and
slightly sharpen) a lemma used in work on the hidden subgroup problem
of quantum computing.

The methods of this paper extend to other algebraic and combinatorial
structures, such as spherical functions of Gelfand pairs, and Bratteli
diagrams. This will be treated in a sequel.

\section{Preliminaries on Markov chain theory} \label{prelim}

This section collects some background on finite Markov chains. Let $X$
be a finite set and $K$ a matrix indexed by $X \times X$ whose rows
sum to 1. Let $\pi$ be a distribution such that $K$ is reversible with
respect to $\pi$; this means that $\pi(x) K(x,y) = \pi(y) K(y,x)$ for
all $x,y$ and implies that $\pi$ is a stationary distribution for the
Markov chain corresponding to $K$. Define the inner product space
$L^2(\pi)$ by letting $<f,g> = \sum_{x \in X} f(x) g(x) \pi(x)$ for
real functions $f,g$. Then when $K$ is considered as an operator on
$L^2(\pi)$ by \[ Kf(x) := \sum_y K(x,y) f(y),\] it is self
adjoint. Hence $K$ has an orthonormal basis of eigenvectors $f_i(x)$
with $Kf_i(x) = \beta_i f_i(x)$, where both $f_i$ and $\beta_i$ are
real. It is easily shown that the eigenvalues satisfy $-1 \leq
\beta_{|X|-1} \leq \cdots \leq \beta_1 \leq \beta_0=1$.

Recall that the total variation distance between probabilities $P,Q$ on $X$ is defined as $||P-Q||_{TV}=\frac{1}{2} \sum_{x \in X} |P(x)-Q(x)|$. It is not hard to see (and
will be used later in proving lower bounds), that
\[ ||P-Q||_{TV} = max_{A \subseteq X} |P(A)-Q(A)| .\] Let $K_x^r$ be the probability measure given by taking $r$ steps from
the starting state $x$. We will be interested in the behavior of
$||K_x^r - \pi||_{TV}$.

The following lemma is well known. Part 1 is the usual method for
computing the power of a diagonalizable matrix. Part 2 upper bounds
$||K_x^r - \pi||_{TV}$ in terms of eigenvalues and eigenvectors and
seems to be remarkably effective in many examples. See \cite{DH} for a
proof of part 2.

\begin{lemma} \label{genbound}
\begin{enumerate}
\item $K^r(x,y) = \sum_{i=0}^{|X|-1} \beta_i^r f_i(x) f_i(y) \pi(y)$ for any $x,y \in X$.
\item  \[ 4 ||K_x^r - \pi||_{TV}^2 \leq \sum_y \frac{|K^r(x,y) - \pi(y)|^2}{\pi(y)} = \sum_{i=1}^{|X|-1} \beta_i^{2r} |f_i(x)|^2 .\] Note that the final sum does not include $i=0$.
\end{enumerate}
\end{lemma}

Finally, let us give a precise definition of the cutoff phenomenon,
taken from \cite{Sal}. Consider a family of finite sets $X_n$, each
equipped with a stationary distribution $\pi_n$, and with another
probability measure $p_n$ that induces a random walk on $X_n$. We say that the cutoff phenomenon holds for the family $(X_n,\pi_n)$ if there exists a sequence $(t_n)$ of positive reals such that 
\begin{enumerate}
\item $\lim_{n \rightarrow \infty} t_n = \infty$;
\item For any $\epsilon \in (0,1)$ and $r_n = [(1+\epsilon)t_n]$, $\lim_{n \rightarrow \infty} ||p_n^{r_n}-\pi_n||_{TV}=0$;
\item For any $\epsilon \in (0,1)$ and $r_n = [(1-\epsilon)t_n]$,
$\lim_{n \rightarrow \infty} ||p_n^{r_n}-\pi_n||_{TV}=1$.
\end{enumerate} The paper \cite{D2} is a nice survey of the cutoff phenomenon.

\section{Preliminaries on group theory} \label{groupsprelim}

This section collects and develops some useful group theoretic
information. Throughout $K$ is the Markov chain on $Irr(G)$ defined
using a representation $\eta$ (not necessarily irreducible) whose
character is real valued.

\begin{lemma} \label{diaggroup} (\cite{F2}, Proposition 2.3) The eigenvalues of $K$ are indexed by conjugacy classes $C$ of $G$:
\begin{enumerate}
\item The eigenvalue parameterized by $C$ is
$\frac{\chi^{\eta}(C)}{d_{\eta}}$.
\item An orthonormal basis of eigenfunctions $f_C$ in $L^2(\pi)$ is
defined by $f_C(\rho) = \frac{|C|^{1/2} \chi^{\rho}(C)}{d_{\rho}}$.
\end{enumerate}
\end{lemma} 

Lemma \ref{decomp} relates the transition probabilities of $K$ to the
decomposition of tensor products. This will be useful in proving the
lower bound for the convergence rate of $K$ when the group is
$GL(n,q)$.

\begin{lemma} \label{decomp} Let $\hat{1}$ denote the trivial representation of $G$. Then \[ K^r(\hat{1},\rho) = \frac{d_{\rho}}{d_{\eta}^r} m_{\rho}(\eta^r),\] where $m_{\rho}(\eta^r)$ denotes the multiplicity of $\rho$ in $\eta^r$. \end{lemma}

\begin{proof} Lemma \ref{diaggroup} and part 1 of Lemma \ref{genbound} imply that \begin{eqnarray*} K^r(\hat{1},\rho) & = & \sum_C  \left( \frac{\chi^{\eta}(C)}{d_{\eta}} \right)^r f_C(\hat{1}) f_C(\rho) \pi(\rho) \\
& = & \sum_C \left( \frac{\chi^{\eta}(C)}{d_{\eta}} \right)^r |C|
\frac{\chi^{\rho}(C)}{d_{\rho}} \frac{d_{\rho}^2}{|G|}\\ & = &
\frac{d_{\rho}}{d_{\eta}^r} \frac{1}{|G|} \sum_g \chi^{\eta}(g)^r
\chi^{\rho}(g). \end{eqnarray*} The result now follows since
$\chi^{\eta}$ is real valued. \end{proof}

Finally, we derive a result which should be useful in many examples
for lower bounding the convergence rate of the Markov chain $K$. It
will be applied in Section \ref{symirrep}.

\begin{prop}
 \label{evolve} Let $C$ be a conjugacy class of $G$ satisfying
 $C=C^{-1}$ and let $f_C$ be as in Lemma \ref{diaggroup}. Let
 $\ee_{K^r}[(f_C)^s]$ denote the expected value of $(f_C)^s$ after r
 steps of the random walk $K$ started at the trivial
 representation. Let $p_{s,C}(T)$ be the probability that the random walk
 on $G$ generated by $C$ and started at the identity is in the
 conjugacy class $T$ after $s$ steps. Then \[ \ee_{K^r}[(f_C)^s] =
 |C|^{s/2} \sum_T p_{s,C}(T) \left( \frac{\chi^{\eta}(T)}{d_{\eta}} \right)^r,\] where the sum is
 over all conjugacy classes $T$ of $G$. \end{prop}

\begin{proof} Let $\hat{1}$ denote the trivial representation of $G$. It follows from Lemma \ref{diaggroup} and part 1 of Lemma \ref{genbound} that \begin{eqnarray*} \ee_{K^r}[ (f_C)^s] & = & \sum_{\rho} K^r(\hat{1},\rho) f_C(\rho)^s \\
& = & \sum_{\rho} \sum_{T} \left( \frac{\chi^{\eta}(T)}{d_{\eta}}
\right)^r f_T(\hat{1}) f_T(\rho) \frac{d_{\rho}^2}{|G|} f_C(\rho)^s
\\ & = & \sum_{\rho} \sum_{T} \left( \frac{\chi^{\eta}(T)}{d_{\eta}}
\right)^r \frac{|T|}{|G|} d_{\rho}^2 \frac{\chi^{\rho}(T)}{d_{\rho}} \left( \frac{|C|^{1/2}
\chi^{\rho}(C)}{d_{\rho}} \right)^s \\ & = & |C|^{s/2} \sum_T \left(
\frac{\chi^{\eta}(T)}{d_{\eta}} \right)^r \frac{|T|}{|G|} \sum_{\rho}
d_{\rho}^2 \frac{\chi^{\rho}(T)}{d_{\rho}} \left(
\frac{\chi^{\rho}(C)}{d_{\rho}} \right)^s. \end{eqnarray*}

  To complete the proof we use the fact that \[ p_{s,C}(T) =
\frac{|T|}{|G|}  \sum_{\rho} d_{\rho}^2
\frac{\chi^{\rho}(T)}{d_{\rho}} \left(
\frac{\chi^{\rho}(C)}{d_{\rho}} \right)^s.\] This is the standard
Fourier analytic expression for $p_{s,C}(T)$; it is explicitly stated
as Exercise 7.67 of \cite{St2} and also follows from Chapter 2 of
\cite{D1}. \end{proof}

\section{The symmetric group} \label{symirrep}

The primary purpose of this section is to obtain sharp convergence
rates and a cutoff phenomenon for random walk on $Irr(S_n)$ when
$\eta$ is the defining representation of $S_n$, whose character on a
permutation is the number of fixed points. Subsection \ref{state}
states and discusses the main result, as well as other interesting
interpretations of the random walk (none of which will be needed for
the proof of the main result). The main result is proved in Subsection
\ref{mainsn}.

\subsection{Main result: statement and discussion} \label{state}

The following theorem is the main result in this paper concerning
random walk on $Irr(S_n)$. For its statement, recall from Section
\ref{prelim} that the total variation distance $||P-Q||_{TV}$ between
two probability distributions $P,Q$ on a finite set $X$ is defined as
$\frac{1}{2} \sum_{x \in X} |P(x)-Q(x)|$. Also note that Theorem
\ref{cutoffsn} proves a cutoff phenomenon (as defined in Section
\ref{prelim}) for the random walk.

\begin{theorem} \label{cutoffsn} Let $G$ be the symmetric group
$S_n$ and let $\pi$ be the Plancherel measure of $G$. Let $\eta$ be
the $n$-dimensional defining representation of $S_n$. Let $K^r$ denote
the distribution of random walk on $Irr(G)$ after $r$ steps, started
from the trivial representation.
\begin{enumerate}
\item If $r=\frac{1}{2}n\log(n)+cn$ with $c \geq 1$ then \[
||K^r-\pi||_{TV} \leq \frac{e^{-2c}}{2}. \]
\item If $r=\frac{1}{2}n\log(n)-cn$ with $0 \leq c \leq \frac{1}{6}
\log(n)$, then there is a universal constant $a$ (independent of $c,n$)
so that \[ ||K^r-\pi||_{TV} \geq 1 -ae^{-4c}.\]
\end{enumerate}
\end{theorem}

It is well known \cite{Sag} that the irreducible representations of
$S_n$ correspond to partitions of $n$, with the trivial representation
corresponding to the one-row partition of size $n$. Thus it is natural
to ask whether the random walk of Theorem \ref{cutoffsn} has a simple
combinatorial description. Proposition \ref{simple} (which was
implicit in \cite{F1}) shows that it does.

\begin{prop} \label{simple}
 A step in the random walk of Theorem \ref{cutoffsn} has the following
 combinatorial description on partitions of $n$. From a partition
 $\lambda$, first one moves to a partition $\mu$ of size $n-1$ which
 can be obtained by removing a corner box from $\lambda$; the chance
 of moving to $\mu$ is $\frac{d_{\mu}}{d_{\lambda}}$. Then one moves
 from $\mu$ to a partition $\rho$ of size $n$ by adding a corner box to
 $\mu$; the chance of moving to $\rho$ is $\frac{d_{\rho}}{n
 d_{\mu}}$. \end{prop}

\begin{proof} By definition, the chance of moving from $\lambda$ to $\rho$ is
 $\frac{d_{\rho}}{nd_{\lambda}}$ multiplied by the multiplicity of
 $\rho$ in $\lambda \otimes \eta$, where $\eta$ is the n-dimensional
 defining representation of $S_n$. Lemma 3.5 of \cite{F1} shows that
 $\lambda \otimes \eta$ is equal to the representation of $S_n$
 obtained by restricting $\lambda$ to $S_{n-1}$ and then inducing it
 to $S_n$. Thus the branching rules for restriction and induction in
 the symmetric group \cite{Sag} imply that the multiplicity of $\rho$
 in $\lambda \otimes \eta$ is the number of $\mu$ of size $n-1$ which
 can be obtained from both of $\lambda,\rho$ by removing some corner
 box (this number is at most 1 if $\lambda \neq \rho$). The result
 follows. \end{proof}

{\bf Remark:} Theorem 3.1 of \cite{F3} gave yet another description of
the random walk on partitions of $n$ corresponding to the walk on
$Irr(S_n)$ in Theorem \ref{cutoffsn}. It proved that the partition
corresponding to a representation chosen from $K^r$ has the same
distribution as the Robinson-Schensted-Knuth (RSK) shape of a
permutation obtained after $r$ iterations of the top to random
shuffle. Hence Theorem \ref{cutoffsn} determines the precise
convergence rate of this RSK shape; note that it is faster than that
of the top to random shuffle, which takes $n\log(n)+cn$ steps to
become random. The reader may wonder why one would care about such
statistics, but the distribution of the RSK shape of a permutation
after various shuffling methods is interesting (see \cite{F6} for more
discussion, including an explanation of why Johansson's work on
discrete orthogonal polynomial ensembles \cite{J} determines the
convergence rate of the RSK shape after riffle shuffles).

To close the discussion of Theorem \ref{cutoffsn}, it is interesting
to compare it with the following classic result on the random
transposition walk on $S_n$.

\begin{theorem} \label{transpos} (\cite{DSh}) Consider
 random walk on the symmetric group $S_n$, where at each step two
symbols are chosen uniformly at random (possibly the same symbol), and
are transposed. If $r=\frac{1}{2}\log(n) + cn$ with $c>0$, then after
$r$ iterations the total variation distance to the uniform
distribution is at most $ae^{-2c}$ for a universal constant
$a$. Moreover, for $c<0$, as $n \rightarrow \infty$, the total
variation distance to the uniform distribution is at least
$(\frac{1}{e}-e^{-e^{-2c}}) + o(1)$. \end{theorem}

Theorems \ref{cutoffsn} and \ref{transpos} establish cutoffs for the
random walks in question, and there is a close parallel between
them. The convergence rates are essentially the same, and whereas
Theorem \ref{transpos} is based on the class of transpositions, which
is the conjugacy class closest to the identity, Theorem \ref{cutoffsn}
is based on the defining representation, which decomposes into the
trivial representation and the representation closest to the trivial
representation. 

\subsection{Proof of Theorem \ref{cutoffsn}} \label{mainsn}

The purpose of this subsection is to prove Theorem
\ref{cutoffsn}. Lemma \ref{fix} is useful for proving the upper bound
in Theorem \ref{cutoffsn}.

\begin{lemma} \label{fix} Suppose that $0 \leq i \leq n-2$. Then the number of permutations on $n$ symbols with exactly $i$ fixed points is at most $\frac{1}{2} \frac{n!}{i!}$.
\end{lemma}

\begin{proof} The number of permutations on $n$ symbols with exactly $i$ fixed points is ${n \choose i} d_{n-i}$ where $d_m$ denotes the number of
 permutations on $m$ symbols with no fixed points. It follows from the
 principle of inclusion-exclusion (see page 67 of \cite{St1}) that
 \[ d_m = m!  \left( 1 - \frac{1}{1!} + \frac{1}{2!} - \frac{1}{3!} +
 \cdots +(-1)^m \frac{1}{m!} \right). \] Thus $d_m \leq \frac{m!}{2}$
 if $m \geq 2$, which implies the result. \end{proof}

Now the upper bound in Theorem \ref{cutoffsn} can be proved.

\begin{proof} (Of part 1 of Theorem \ref{cutoffsn}) From Lemma \ref{diaggroup} and
 part 2 of Lemma \ref{genbound} it follows that \[ 4
 ||K^r-\pi||_{TV}^2 \leq \sum_{i=0}^{n-2} |\{g \in S_n : fp(g)=i \}|
 \left( \frac{i}{n} \right)^{2r} \] where $fp(g)$ is the number of
 fixed points of $g$. Note that the sum ends at $n-2$ since no
 permutation can exactly have $n-1$ fixed points. Applying Lemma
 \ref{fix} and then letting $j=n-i$, one concludes that \begin{eqnarray*} 4
 ||K^r-\pi||_{TV}^2 & \leq & \frac{1}{2} \sum_{i=0}^{n-2} \frac{n!}{i!}
 \left( \frac{i}{n} \right)^{2r}\\
& = & \frac{1}{2} \sum_{j=2}^n \frac{n!}{(n-j)!}  \left(1-\frac{j}{n}
 \right)^{2r}\\ & = & \frac{1}{2} \sum_{j=2}^n \frac{n!}{(n-j)!} e^{2r
 \cdot \log(1-j/n)}\\ & \leq & \frac{1}{2} \sum_{j=2}^n
 \frac{n!}{(n-j)!} e^{-2rj/n}. \end{eqnarray*}
 Recalling that $r=\frac{1}{2}n\log(n)+cn$, the bound becomes
 \[ \frac{1}{2} \sum_{j=2}^n \frac{n!}{(n-j)!
 n^{j}} e^{-2cj} \leq \frac{1}{2} \sum_{j=2}^n e^{-2cj} =
 \frac{e^{-4c}}{2(1-e^{-2c})}.\] Dividing by 4 and taking square
 roots, the result follows since $c \geq 1$. \end{proof}

The idea for proving the lower bound of Theorem \ref{cutoffsn} will be
to find a random variable on $Irr(G)$ which is far from its
distribution under Plancherel measure if fewer than
$\frac{1}{2}n\log(n)$ steps have been taken. The random variable to be
used is precisely the $f_C$ from Lemma \ref{diaggroup}, where $C$ is
the conjugacy class of transpositions. This approach is dual to, and
motivated by, the approach on page 44 of \cite{D1} for proving the
lower bound in Theorem \ref{transpos}.

\begin{proof} (Of part 2 of Theorem \ref{cutoffsn}). We apply Chebyshev's
 inequality. Let $C$ be the conjugacy class of transpositions, and let
 $f_C$ be as in Lemma \ref{diaggroup}. For $\alpha>0$ to be specified
 later, let $A$ be the event that $f_C \leq \alpha$. The
 orthogonality relations for the irreducible characters of the
 symmetric group imply that under Plancherel measure $\pi$, the
 random variable $f_C$ has mean 0 and variance 1. Hence $\pi(A) \geq
 1-\frac{1}{\alpha^2}$.

On the other hand, it follows from Proposition \ref{evolve} that
\[ \ee_{K^r}[f_C] = \sqrt{ {n \choose 2}} \left( 1-\frac{2}{n}
\right)^r.\] Letting $r=\frac{1}{2}n\log(n)-cn$ and using the Taylor
expansion for $\log(1-x)$ gives that
\begin{eqnarray*} \ee_{K^r}[f_C] & = & \sqrt{ {n \choose 2}}
\exp \left( (\frac{1}{2}n\log(n)-cn) \cdot \log(1-2/n) \right) \\
& = & \sqrt{{n \choose 2}} \exp \left(
-\log(n)+2c+O(\frac{\log(n)}{n})+ O(\frac{c}{n}) \right)\\
 & \geq & \frac{1}{2} \exp \left( 2c+O(\frac{\log(n)}{n}) + O(\frac{c}{n})
\right). \end{eqnarray*} This is large when $c$ is large.

Proposition \ref{evolve} gives that \[ \ee_{K^r}[(f_C)^2] = 1+ {n-2
\choose 2} (1-\frac{4}{n})^r+(2n-4)(1-\frac{3}{n})^r.\] Thus the variance of $f_C$ under $K^r$ is
\begin{eqnarray*} & & 1+ {n-2 \choose 2}
(1-\frac{4}{n})^r+(2n-4)(1-\frac{3}{n})^r - {n \choose 2}
(1-\frac{2}{n})^{2r}\\ & = & 1+ {n-2 \choose 2} \exp \left(
(\frac{1}{2}n\log(n)-cn)(-\frac{4}{n}+O(\frac{1}{n^2}) )\right)\\ & & +
(2n-4) \exp \left(
(\frac{1}{2}n\log(n)-cn)(-\frac{3}{n}+O(\frac{1}{n^2}) ) \right)\\ & &
- {n \choose 2} \exp \left(
(n\log(n)-2cn)(-\frac{2}{n}+O(\frac{1}{n^2})) \right) \\ & = & 1 +
\frac{1}{2} \exp \left( 4c+O(\frac{\log(n)}{n}) + O(\frac{c}{n})
\right)\\ & & + \frac{2}{\sqrt{n}} \exp \left( 3c+O(\frac{\log(n)}{n}) +
O(\frac{c}{n}) \right)\\ & & - \frac{1}{2} \exp \left(
4c+O(\frac{\log(n)}{n}) + O(\frac{c}{n}) \right)\\ & = & 1 +
\frac{e^{4c}}{2} \left( O(\frac{\log(n)}{n}) + O(\frac{c}{n}) \right) +
\frac{2}{\sqrt{n}} \exp \left( 3c+O(\frac{\log(n)}{n}) + O(\frac{c}{n})
\right). \end{eqnarray*} Since $0 \leq c \leq \frac{1}{6} \log(n)$, the
variance is bounded by a universal constant. Let
$\alpha=\frac{e^{2c}}{4}$. Then Chebyshev's inequality gives that
$K^r(A) \leq \frac{b}{e^{4c}}$ for a universal constant $b$. Thus \[
||K^r - \pi||_{TV} \geq |\pi(A) - K^r(A)| \geq 1 -
\frac{1}{\alpha^2}-\frac{b}{e^{4c}}, \] which completes the
proof. \end{proof}

\section{The general linear group} \label{glirrep}

This section studies random walk on $Irr(GL(n,q))$ in the case that
$\eta$ is the representation of $GL(n,q)$ whose character is
$q^{d(g)}$, where $d(g)$ is the dimension of the fixed space of $g$.
Subsection \ref{maingl} states and discusses the main
result. Subsection \ref{upp} proves the upper bound in Theorem
\ref{cutoffgl} and Subsection \ref{low} proves the lower bound.

\subsection{Main result: statement and discussion} \label{maingl}

The following theorem is the main result in the paper concerning
random walk on $Irr(GL(n,q))$.

\begin{theorem} \label{cutoffgl} Let
 $G$ be the finite general linear group $GL(n,q)$ and let $\pi$ be the
Plancherel measure of $G$. Let $\eta$ be the representation of $G$
whose character is $q^{d(g)}$, where $d(g)$ is the dimension of the
fixed space of $g$. Let $K^r$ denote the distribution of random walk
on $Irr(G)$ after $r$ steps started from the trivial representation.
\begin{enumerate}
\item If $r=n+c$ with $c > 0$, then $||K^r-\pi||_{TV} \leq
\frac{1}{2q^c}$.
\item If $r=n-c$ with $c > 0$, then $||K^r-\pi||_{TV} \geq
1-\frac{a}{q^c}$ where $a$ is a universal constant (independent of
$n,q,c$).
\end{enumerate}
\end{theorem}

It is interesting to compare this result with the following result of
Hildebrand on the random transvection walk.

\begin{theorem} \label{transvec} (\cite{H}) Consider random
walk on the finite special linear group $SL(n,q)$, where at each step
one multiplies by a random transvection (i.e. an element of $SL(n,q)$
which is not the identity but fixes all points in a hyperplane). There
are positive constants $a,b$ such that for sufficiently large $n$ and
all $c>0$, the total variation distance to the uniform distribution
after $n+c$ steps is at most $ae^{-bc}$. Moreover, given $\epsilon>0$,
there exists $c>0$ so that the total variation distance is at least
$1-\epsilon$ after $n-c$ iterations for sufficiently large
$n$. \end{theorem}

In both theorems there is a cutoff around $n$ steps. Also, whereas
Theorem \ref{transvec} is based on the class of transvections, which
is the unipotent class closest to the identity, Theorem \ref{cutoffgl}
is based on a representation which for $q=2$ decomposes into the
trivial representation and the unipotent representation closest to the
trivial representation (for $q>2$ the decomposition involves a few
more pieces).

\subsection{Upper bound on convergence rate} \label{upp}

To prove the upper bound of Theorem \ref{cutoffgl}, the following
lemma will be helpful. Recall that $d(g)$ is the dimension of the
fixed space of $g$. We also use the notation that $(1/q)_r = (1-1/q)
\cdots (1-1/q^r)$.

\begin{lemma} \label{fixgl} Suppose that $0 \leq i \leq n$. Then the number of elements in $GL(n,q)$ with $d(g)=i$ is at most $q^{n^2-i^2}$.
\end{lemma}

\begin{proof} It
 is known (going back at least to \cite{RS}) that the number of
 elements in $GL(n,q)$ with $d(g)=i$ is exactly \[
 \frac{|GL(n,q)|}{|GL(i,q)|} \sum_{j=0}^{n-i} \frac{(-1)^j q^{{j
 \choose 2}}}{q^{ij} |GL(j,q)|} .\] Since $|GL(n,q)| =
 q^{n^2}(1/q)_n$, clearly $ \frac{|GL(n,q)|}{|GL(i,q)|} \leq
 q^{n^2-i^2}$. Also \[ \sum_{j=0}^{n-i} \frac{(-1)^j q^{{j \choose
 2}}}{q^{ij} |GL(j,q)|} = \sum_{j=0}^{n-i} \frac{(-1)^j}{q^{ij}
 (q^j-1) \cdots (q-1)} \] is an alternating sum of decreasing terms
 the first of which is 1, so the sum is at most 1. This proves the
 lemma. \end{proof}

\begin{proof} (Of part 1 of Theorem \ref{cutoffgl}) By Lemma \ref{diaggroup} and part 2 of Lemma \ref{genbound}, it follows that \begin{eqnarray*} 4 ||K^r-\pi||_{TV}^2 & \leq & \sum_{i=0}^{n-1} |\{g \in GL(n,q): d(g)=i \}| \left( \frac{1}{q^{n-i}} \right)^{2r}\\
& = & \sum_{i=1}^n |\{g \in GL(n,q): d(g)=n-i \}| \left(
\frac{1}{q^{i}} \right)^{2r}, \end{eqnarray*} where $d(g)$ is the
dimension of the fixed space of $g$. By Lemma \ref{fixgl} this is at
most $\sum_{i=1}^n \frac{q^{2ni}}{q^{i^2+2ir}}$. Since $r=n+c$, this
is equal to \[ \sum_{i=1}^n \frac{1}{q^{i^2+2ci}} \leq q^{-2c}
\sum_{i=1}^n \frac{1}{q^{i^2}} \leq \frac{1}{q(1-1/q)} q^{-2c}.\]
Since $q \geq 2$, this is at most $q^{-2c}$. Dividing by 4 and taking
square roots completes the proof. \end{proof}

\subsection{Lower bound on convergence rate} \label{low}

To prove the lower bound of Theorem \ref{cutoffgl}, we will need to
know about representation theory of $GL(n,q)$ and about the
decomposition of $\eta^k$ into irreducibles, where $1 \leq k \leq
n$. In fact the paper \cite{GK} studies this decomposition.

To begin we recall some facts about $Irr(GL(n,q))$. A full treatment
of the subject with proofs appears in \cite{M}, \cite{Z} but we will
adhere to the notation of \cite{GK} instead. As usual a partition
$\lambda=(\lambda_1,\cdots,\lambda_m)$ is identified with its
geometric image $\{(i,j): 1 \leq i \leq m, 1 \leq j \leq \lambda_i \}$
and $|\lambda|=\lambda_1+\cdots+\lambda_m$ is the total number of
boxes. Let $\yy$ denote the set of all partitions, including the
empty partition of size 0.

Given an integer $1 \leq k < n$ and two characters $\chi_1,\chi_2$ of
the groups $GL(k,q)$ and $GL(n-k,q)$, their parabolic induction
$\chi_1 \circ \chi_2$ is the character of $GL(n,q)$ induced from the
parabolic subgroup of elements of the form \[ P = \left \{ \left( \begin{array} {cc} g_1 & * \\
0 & g_2 \end{array} \right) : g_1 \in GL(k,q), g_2 \in GL(n-k,q) \right \} \] by the function $\chi_1(g_1) \chi_2(g_2)$.

A character is called cuspidal if it is not a component of any
parabolic induction. Let ${\it C}_d$ denote the finite set of cuspidal
characters of $GL(d,q)$ and let ${\it C} = \cup_{d \geq 1} {\it
C}_d$. The unit character of $GL(1,q)$ plays an important role and
will be denoted $e$; it is one of the $q-1$ elements of ${\it
C}_1$. Given a family ${\it \phi}: {\it C} \mapsto \yy$ with finitely
many non-empty partitions ${\it \phi}(c)$, its degree $||\phi||$ is
defined as $\sum_{d \geq 1} \sum_{c \in {\it C}_d} d \cdot |{\it
\phi}(c)|$. We also write $||c||=d$ if $d$ is the unique number so
that $c \in {\it C}_d$. A fundamental result is that the irreducible
representations of $GL(n,q)$ are in bijection with the families of
partitions of degree $n$, so we also let $\phi$ denote the
corresponding representation. The partition ${\mathcal \phi}(e)$ will
be referred to as the unipotent part of ${\mathcal \phi}$.

It will be helpful to know that the dimension of the irreducible
representation of $GL(n,q)$ corresponding to the family ${\mathcal
\phi}$ is \[ (q^n-1) \cdots (q-1) \prod_{d \geq 1} \prod_{c \in
{\mathcal C_d}} \frac{q^{d \cdot n(\phi(c))}}{\prod_{b \in \phi(c)}
(q^{d \cdot h(b)}-1)}, \] where $h(b)$ is the hooklength
$\lambda_i+\lambda_j'-i-j+1$ of a box $b=(i,j)$ and $n(\lambda)=\sum_i
(i-1) \lambda_i$. In what follows, we also use that notation that
$\lambda'$ is the transpose of $\lambda$, obtained by switching the
rows and columns of $\lambda$ (i.e. $\lambda'_i = |\{j: \lambda_j \geq
i\}|$).

\begin{prop} \label{firstrow} Suppose that $K^r(\hat{1},{\mathcal \phi}) > 0$. Then $\phi(e)_1 \geq n-r$.
\end{prop}

\begin{proof} If $r=0$ or $r>n$ the result is trivially true (and not useful to us). So suppose that $1 \leq r \leq n$. Lemma \ref{decomp} implies that $K^r(\hat{1},{\mathcal \phi}) > 0$ if and only if ${\mathcal \phi}$ occurs as a component of $\eta^r$. Proposition 5 and Theorem 7 of \cite{GK} imply that if $1 \leq r \leq n$ and  ${\mathcal \phi}$ occurs as a component of $\eta^r$, then $\phi(e)_1 \geq n-r$. \end{proof}  

For $c>0$, define $A$ as the event that a representation has the first
row of its unipotent part of size at least $c$. Proposition
\ref{firstrow} showed that if $r=n-c$, then $K^r(A)=1$. The next goal
will be to upper bound $\pi(A)$ where $\pi$ is the Plancherel measure
of $GL(n,q)$. First we upper bound the Plancherel probability that the
unipotent part of a random representation is $\lambda$.

\begin{lemma} \label{planbound} Let $\pi$ be the Plancherel measure of $GL(n,q)$. Then for any $\lambda$, \[ \pi(\phi(e)=\lambda) \leq \frac{1}{q^{\sum_i (\lambda_i)^2} \prod_{b \in \lambda} (1-1/q^{h(b)})^2}.\]
\end{lemma}

\begin{proof} By the definition of Plancherel measure and the formula for $d_{\phi}$ one has that \begin{eqnarray*} \pi(\phi(e)=\lambda) & = & \sum_{||\phi||=n \atop \phi(e)=\lambda} \frac{d_{\phi}^2}{|GL(n,q)|}\\
& = & \frac{(q^n-1) \cdots (q-1)}{q^{{n \choose 2}}} \sum_{||\phi||=n
\atop \phi(e)=\lambda} \prod_{d \geq 1} \prod_{c \in C_d} \frac{q^{2d
\cdot n(\phi(c))}}{ \prod_{b \in \phi(c)} (q^{d \cdot h(b)}-1)^2}\\ &
= & \frac{q^{n+2n(\lambda)} (1/q)_n}{\prod_{b \in \lambda}
(q^{h(b)}-1)^2} \sum_{||\phi||={n-|\lambda|} \atop \phi(e)=\emptyset}
\prod_{d \geq 1} \prod_{c \in C_d} \frac{q^{2d \cdot n(\phi(c))}}{
\prod_{b \in \phi(c)} (q^{d \cdot h(b)}-1)^2}\\ & \leq &
\frac{q^{n+2n(\lambda)} (1/q)_n}{\prod_{b \in \lambda} (q^{h(b)}-1)^2}
\sum_{||\phi||={n-|\lambda|}} \prod_{d \geq 1} \prod_{c \in C_d}
\frac{q^{2d \cdot n(\phi(c))}}{ \prod_{b \in \phi(c)} (q^{d \cdot
h(b)}-1)^2}\\ & = & \frac{q^{n+2n(\lambda)} (1/q)_n}{\prod_{b \in
\lambda} (q^{h(b)}-1)^2} \sum_{ ||\phi||={n-|\lambda|}}
\frac{d_{\phi}^2}{(q^{n-|\lambda|}-1)^2 \cdots
(q-1)^2}. \end{eqnarray*} Since the sum of the squares of the
dimensions of the irreducible representations of a finite group is
equal to the order of the group, this is \begin{eqnarray*} & & \frac{q^{n+2n(\lambda)}
(1/q)_n}{\prod_{b \in \lambda} (q^{h(b)}-1)^2}
\frac{|GL(n-|\lambda|,q)|}{(q^{n-|\lambda|}-1)^2 \cdots (q-1)^2}\\ & =
& \frac{q^{2n(\lambda)+|\lambda|} (1/q)_n}{ (1/q)_{n-|\lambda|}
\prod_{b \in \lambda} (q^{h(b)}-1)^2}\\ & \leq &
\frac{q^{2n(\lambda)+|\lambda|}}{ \prod_{b \in \lambda}
(q^{h(b)}-1)^2}\\ & = & \frac{1}{q^{|\lambda|+2n(\lambda')} \prod_{b
\in \lambda} (1-1/q^{h(b)})^2}\\ & = & \frac{1}{q^{\sum_i
(\lambda_i)^2} \prod_{b \in \lambda} (1-1/q^{h(b)})^2}
. \end{eqnarray*} The second to last equation used the identity \[
\sum_{b \in \lambda} h(b) = n(\lambda) + n(\lambda') + |\lambda| \] on
page 11 of \cite{M}.
\end{proof}

\begin{prop} \label{planunip} There is a universal constant $a$ (independent of $n,q,c$) so that $\pi(|\phi(e)| \geq c) \leq \frac{a}{q^c}$.
\end{prop}

\begin{proof} By Lemma \ref{planbound}, \[ \sum_{|\lambda|=m} \pi(\phi(e)=\lambda) \leq  \sum_{|\lambda|=m} \frac{1}{q^{\sum_i (\lambda_i)^2} \prod_{b
\in \lambda} (1-1/q^{h(b)})^2}.\] Noting that this sum is invariant under transposing $\lambda$, Proposition 4.2 of
\cite{F5} implies that it is equal to the coefficient of $u^m$ in
$\prod_{i=1}^{\infty} \prod_{j=0}^{\infty} \left(
\frac{1}{1-\frac{u}{q^{i+j}}} \right)$. Lemma 4.4 of \cite{F5} shows
that since $q \geq 2$, this coefficient is at most
$\frac{1}{(q^m-1)(1-1/q)^6}$. Thus \[ \pi(|\phi(e)| \geq c) = \sum_{m
\geq c} \pi(|\phi(e)|=m) \leq \frac{1}{(1-1/q)^6} \sum_{m \geq c}
\frac{1}{q^m-1},\] which implies the result. \end{proof}

The lower bound of Theorem \ref{cutoffgl} can now be proved.

\begin{proof} (Of part 2 of Theorem \ref{cutoffgl}). Fix $c>0$ and define $A$ as the event that a representation has the first row of its unipotent part of size at least $c$. Defining $r=n-c$, Proposition \ref{firstrow} showed that
$K^r(A)=1$. Clearly $\pi(A) \leq \pi(|\phi(e)| \geq c)$. By
Proposition \ref{planunip}, this is at most $\frac{a}{q^c}$ for a
universal constant $a$. Since $||K^r-\pi||_{TV} \geq
|\pp_{K^r}(A)-\pp_{\pi}(A)|$, the result follows. \end{proof}

\section{Asymptotic description of Plancherel measure of $GL(n,q)$}
 \label{plancherel}

Given the numerous papers on asymptotics of Plancherel measure of the
symmetric groups (see \cite{J},\cite{O},\cite{BOO} and references
therein), it is natural to study asymptotics of Plancherel measure for
other towers of finite groups. Aside from the paper \cite{VK}, which
is only tangentially related and contains no proofs, we are aware of
no results on this question. In this section an elegant asymptotic
description of the Plancherel measure of $GL(n,q)$ is obtained, when
$q$ is fixed and $n \rightarrow \infty$: it is proved that this
limiting measure factors into independent pieces, the distributions of
which can be explicitly described. Then a connection with the proof of
the convergence rate lower bound of Theorem \ref{cutoffgl} is noted.

Fix $c \in {\mathcal C}$ and let ${\mathcal \phi}$ be a representation
chosen from the Plancherel measure $\pi$ of $GL(n,q)$. Then the
partition $\phi(c)$ is a random partition (of size at most $n$). It
is natural to study the distribution of $\phi(c)$ when $c$ is fixed
and $n \rightarrow \infty$. Theorem \ref{asymp} will show that the
random partitions $\{ \phi(c): c \in {\mathcal C} \}$ are independent
in the $n \rightarrow \infty$ limit, and will determine the
distribution of each of them.

First, we define a ``cycle index'' $\hat{Z}_{GL(n,q)}$ for irreducible
representations. For $n \geq 1$, let \[ \hat{Z}_{GL(n,q)} =
\sum_{\phi: ||\phi||=n} \pi(\phi) \prod_{c \in {\mathcal C}: |\phi(c)|
\neq 0} x_{c,\phi(c)}.\] Here the $x_{c,\phi(c)}$ are variables
corresponding to pairs of elements of ${\mathcal C}$ and partitions.

\begin{lemma} \label{cycindex}
\[ 1 + \sum_{n=1}^{\infty} \hat{Z}_{GL(n,q)} \frac{u^n}{(1/q)_n} = \prod_{d \geq 1} \prod_{c \in {\mathcal C_d}} \left[ 1 + \sum_{|\lambda| \geq 1} \frac{ x_{c,\lambda} \cdot u^{d |\lambda|}}{q^{d \sum_i (\lambda_i)^2} \prod_{b \in \lambda} (1-1/q^{d \cdot h(b)})^2} \right]\]
\end{lemma}

\begin{proof} From the formula for the dimension of an element of $Irr(GL(n,q))$ given
 in Subsection \ref{low}, it follows by comparing the coefficients of
 products of the x variables on both sides that \[  1 + \sum_{n=1}^{\infty} \hat{Z}_{GL(n,q)} \frac{u^n}{q^n (1/q)_n} = \prod_{d \geq 1}
 \prod_{c \in {\C_d}} \left[ 1 + \sum_{|\lambda| \geq 1}
 \frac{x_{c,\lambda} \cdot u^{d|\lambda|} q^{2d \cdot n(\lambda)} }{ \prod_{b
 \in \lambda} (q^{d \cdot h(b)}-1)^2} \right].\] Replacing $u$ by $uq$ gives that \[  1 + \sum_{n=1}^{\infty} \hat{Z}_{GL(n,q)} \frac{u^n}{(1/q)_n} = \prod_{d \geq 1}
 \prod_{c \in {\C_d}} \left[ 1 + \sum_{|\lambda| \geq 1}
 \frac{x_{c,\lambda} \cdot u^{d|\lambda|} q^{d(|\lambda|+2n(\lambda))} }{ \prod_{b
 \in \lambda} (q^{d \cdot h(b)}-1)^2} \right].\] Arguing as in the
 last two lines of the proof of Lemma \ref{planbound} proves the
 result. \end{proof}

To state the main result of this subsection, we define, for $q>1$ and
$0<u<q$, a probability measure $S_{u,q}$ on $\yy$ (the set of all
partitions of all natural numbers). This is defined by the formula \[
S_{u,q}(\lambda) = \prod_{i=1}^{\infty} \prod_{j=0}^{\infty} \left(1 -
\frac{u}{q^{i+j}} \right) \cdot \frac{u^{|\lambda|}}{q^{\sum_i (\lambda_i)^2}
\prod_{b \in \lambda} (1-1/q^{h(b)})^2}.\] This measure, and some of
its properties, are discussed in \cite{F5} (note that there $\lambda$
is replaced by its transpose).

The following simple lemma is useful.

\begin{lemma} \label{taylor} If a function $f(u)$ has a Taylor series around 0 which converges at $u=1$, then the $n \rightarrow \infty$ limit of the coefficient of $u^n$ in $\frac{f(u)}{1-u}$ is equal to $f(1)$.
\end{lemma}

\begin{proof} Write the Taylor series $f(u)=\sum_{n=0}^{\infty} a_n u^n$. Then observe that the coefficient of $u^n$ in $\frac{f(u)}{1-u}$ is equal to $\sum_{i=0}^n a_i$. \end{proof}

Now the main theorem of this subsection can be proved.

\begin{theorem} \label{asymp}
\begin{enumerate}
\item Fix $u$ with $0<u<1$. Then choose a random natural number $N$
with $\pp(N=n) = \prod_{m=0}^{\infty} \left( 1 - \frac{u}{q^m} \right)
\cdot \frac{u^n}{(1/q)_n}$. Choose $\phi$ from the Plancherel measure
of $GL(N,q)$. Then as $c \in {\mathcal C}$ varies, the random
partitions $\phi(c)$ are independent with $\phi(c)$ distributed
according to the measure $S_{u^{||c||},q^{||c||}}$.

\item Choose $\phi$ from the Plancherel measure $\pi_n$ of
$GL(n,q)$. Then as $n \rightarrow \infty$, the random partitions
$\phi(c)$ converge to independent random variables, with $\phi(c)$
distributed according to the measure $S_{1,q^{||c||}}$.
\end{enumerate}
\end{theorem}

\begin{proof} Setting all of the variables $x_{c,\lambda}$ equal to 1 in Lemma
 \ref{cycindex}, the left hand side becomes $\sum_{n \geq 0}
 \frac{u^n}{(1/q)_n}$, which by an identity of Euler (Corollary 2.2 of
 \cite{An}) is equal to $\prod_{m=0}^{\infty} (1-u/q^m)^{-1}$. Since
 $S_{u,q}$ is a probability measure, the right hand side becomes $
 \prod_{d \geq 1} \prod_{c \in {\mathcal C_d}} \prod_{i=1}^{\infty}
 \prod_{j=0}^{\infty} \left( 1-\frac{u^{d}}{q^{d(i+j)}}
 \right)^{-1}$. Taking reciprocals of this equation and multiplying
 by the statement of Lemma \ref{cycindex} gives the
 equation \begin{eqnarray*} & & \prod_{m=0}^{\infty} (1-u/q^m) +
 \sum_{n=1}^{\infty} \hat{Z}_{GL(n,q)} \prod_{m=0}^{\infty} (1-u/q^m)
 \cdot \frac{u^n}{(1/q)_n} \\ & = & \prod_{d \geq 1} \prod_{c \in
 {\mathcal C_d}} \left( S_{u^d,q^d}(\emptyset) + \sum_{|\lambda| \geq
 1} S_{u^d,q^d}(\lambda) x_{c,\lambda} \right). \end{eqnarray*} This
 proves the first assertion.

 For the second assertion, divide both sides of the previous equation
 by $\prod_{m=1}^{\infty} (1-u/q^m)$, giving that \begin{eqnarray*} &
 & (1-u) \left( 1 + \sum_{n=1}^{\infty} \hat{Z}_{GL(n,q)}
 \frac{u^n}{(1/q)_n} \right)\\ & = & \prod_{m=1}^{\infty}
 (1-u/q^m)^{-1} \prod_{d \geq 1} \prod_{c \in {\mathcal C_d}} \left(
 S_{u^d,q^d}(\emptyset) + \sum_{|\lambda| \geq 1} S_{u^d,q^d}(\lambda)
 x_{c,\lambda} \right). \end{eqnarray*} Thus for any $c_1,\cdots,c_t
 \in {\mathcal C}$ and any $\lambda_1,\cdots,\lambda_t \in \yy$, it
 follows that \[ \lim_{n \rightarrow \infty} \
 \pi_n(\phi(c_1)=\lambda_1,\cdots, \phi(c_t)=\lambda_t) \] is equal to
 the limit as $n \rightarrow \infty$ of \[ (1/q)_n \cdot [u^n]
 \frac{1}{1-u} \prod_{m=1}^{\infty} (1-u/q^m)^{-1} \prod_{i=1}^t
 S_{u^{||c_i||},q^{||c_i||}}(\lambda_i),\] where $[u^n] f(u)$ denotes
 the coefficient of $u^n$ in $f(u)$. By Lemma \ref{taylor} this limit
 exists and is $\prod_{i=1}^t S_{1,q^{||c_i||}}(\lambda_i)$, as
 desired. \end{proof}

To close this section, note that part 2 of Theorem \ref{asymp} gives
an intuitive explanation of why Proposition \ref{planunip} should be
true. Indeed, $\phi(e)$ converges to a partition chosen from $S_{1,q}$
as $n \rightarrow \infty$, and a partition chosen from $S_{1,q}$ has
finite size. Thus when $c$ is big, the probability that $|\phi(e)|
\geq c$ should be small.	

\section{Connection to quantum computing} \label{quantum}

As is explained in Chapter 5 of the text \cite{NC}, many of the
problems in which a quantum computer outperforms its classical
counterpart, such as factoring and the discrete-log problem, can be
described in terms of the following hidden subgroup problem. Let $G$
be a finite group and $H$ a subgroup. Given a function $f$ from $G$ to
a finite set that is constant on left cosets $gH$ of $H$ and takes
different values for different cosets, the hidden subgroup problem is
to determine a set of generators of $H$ and the decision version of
the problem is to determine whether there is a non-identity hidden
subgroup or not.

One approach to these problems uses the ``weak standard method'' of
quantum Fourier sampling, which is described in \cite{KS} and more
fully in \cite{HRT}. When applying this to a subgroup $H$, one obtains
a probability measure $P_H$ on $Irr(G)$, which chooses a
representation $\rho \in Irr(G)$ with probability \[ P_H(\rho) =
\frac{d_{\rho}}{|G|} \sum_{h \in H} \chi^{\rho}(h).\] Then a subgroup
$H$ can be distinguished efficiently from the trivial subgroup $\{e\}$
if and only if the total variation distance $||P_H-P_{\{e\}}||_{TV}$
is larger than $(\log |G|)^{-c}$ for some constant $c$.

The following upper bound on the total variation distance
$||P_H-P_{\{e\}}||_{TV}$ was used by Kempe and Shalev in their work on
the weak standard method, and improved earlier results in the
literature.

\begin{prop} \label{kemp} (\cite{KS})
 Let $C_1,\cdots,C_k$ denote the non-identity conjugacy classes of
 $G$. Then \[ ||P_H-P_{\{e\}}||_{TV} \leq \frac{1}{2} \sum_{i=1}^k |C_i \cap H|
 |C_i|^{-1/2}.\] \end{prop}

We use the perspective of random walk on $Irr(G)$ to prove Proposition
\ref{better}, which is a slight sharpening of Proposition
\ref{kemp}. (To see that it is sharper, take squares). Note also that
the random walk viewpoint ``explains'' the appearance of the
quantities $|C_i \cap H| |C_i|^{-1}$ in the total variation upper
bounds of Propositions \ref{kemp} and \ref{better}: they are simply
the eigenvalues of random walk on $Irr(G)$.

\begin{prop} \label{better} Let $C_1,\cdots,C_k$ denote the non-identity conjugacy classes of
 $G$. Then \[ ||P_H-P_{\{e\}}||_{TV} \leq \frac{1}{2} \left[
 \sum_{i=1}^k |C_i \cap H|^2 |C_i|^{-1} \right]^{1/2}.\] \end{prop}

\begin{proof} Note that $P_{\{e\}}$ is simply the
Plancherel measure $\pi$ on $Irr(G)$. Next define $\eta$ to be the
induced representation $Ind_H^G(\hat{1})$, where $\hat{1}$ denotes the
trivial representation of $H$. Then the character of $\eta$ is real
valued, since its value on $g$ is the number of left cosets of $H$
fixed by $g$. The probability that random walk on $Irr(G)$ defined by
$\eta$ and started at the trivial representation is at $\rho$ after
one step is equal to $\frac{d_{\rho}}{d_{\eta}}
m_{\rho}(Ind_H^G(\hat{1}))$, where $m_{\rho}(Ind_H^G(\hat{1}))$ is the
multiplicity of $\rho$ in $Ind_H^G(\hat{1})$. Frobenius reciprocity
gives that $m_{\rho}(Ind_H^G(\hat{1}))$ is equal to $\frac{1}{|H|}
\sum_{h \in H} \chi^{\rho}(h)$. Since $d_{\eta}=\frac{|G|}{|H|}$, it
follows that the chance of being at $\rho$ after 1 step of the random
walk on $Irr(G)$ is precisely equal to $P_H(\rho)$. Summarizing, \[
||P_H-P_{\{e\}}||_{TV} = ||K^1 - \pi||_{TV}, \] where $K^1$ is the
distribution on $Irr(G)$ after 1 step started from the trivial
representation.

Next, apply part 2 of Lemma \ref{genbound} with $r=1$ and Lemma
\ref{diaggroup} to conclude that \[ ||K^1 - \pi||_{TV} \leq
\frac{1}{2} \left[ \sum_{i=1}^k \left(
\frac{\chi^{\eta}(C_i)}{d_{\eta}} \right)^2 |C_i| \right]^{1/2} .\]
From the formula for induced characters (page 47 of \cite{Sag}), it
follows that $\frac{\chi^{\eta}(C_i)}{d_{\eta}} = \frac{|C_i \cap
H|}{|C_i|}$, which completes the proof. \end{proof}

\end{document}